\documentclass[12pt]{article}
\usepackage{amssymb}
\usepackage{amsmath}
\usepackage{latexsym}
\newtheorem{theorem}{Theorem}[section]

\newcommand {\ao} { { {} \atop 1 } }
\newcommand {\ad} { { {} \atop 2 } }

\def \dsps {\displaystyle}
\def \lngra {\longrightarrow}

\def \hgo {\widehat{g^1} }
\def \hgdv {\widehat{g^2} }

\def \Xsi {{\cal X}_i}
\def \vchio { V_{\chi_\ao} }
\def \vchidv { V_{\chi_\ad} }

\def \chio { \chi_\ao }
\def \chidv { \chi_\ad }

\def \hg  { {\widehat g} }

\def \hojch { H_{2,J(p)} (S, V_\chi  \otimes \Delta) }

\def \hso { H_{2,J_1(p)} (S_1, \vchio \otimes \Delta_1) }
\def \hsdv { H_{2,J_2(p)} (S_2, V_{\chi_\ad} \otimes \Delta_2) }

\def \hsjo { H_{2,J_{1,0},...,J_{1,k-1} } (S_1, \vchio \otimes \Delta_1) }
\def \hsjdv { H_{2, J_{2,0},...,J_{2,k-1} }
 (S_2, V_{\chi_\ad} \otimes \Delta_2) }

\def \hf { {\widehat f   }}
\def \hg { {\widehat g   }}

\def \htwoch {H_2(S, V_\chi \otimes \Delta)}

\def \bundo{ V_{\chi_\ao}\otimes \Delta_1 }
\def \bundv { V_{\chi_\ad } \otimes \Delta_2 }

\def \hjb { H_{ 2, J_{0},...J_{k-1} } (S, V_\chi\otimes \Delta ) }
\def \hjo { H_{2,J_{1,0}, ..., J_{1,k-1} }(S_1, V_{\chi_\ao}
 \otimes \Delta_1 )}
\def \hjdv  {H_{2,J_{2,0} ,..., J_{2,k-1} }(S, V_{\chi_\ad} \otimes \Delta_2 )}

\def \Dlt {\Delta}

\def \ttilp { T_{ {\widetilde p} } }

\def \hhdv  { {\widehat h}^2 }
\def \hho   { {\widehat h}^1 }

\def \lmpt {\longrightarrow}

\def \tilp {{\widetilde p}}
\def \piod {\pi_1}

\def \hfo { {\widehat f^1 } }
\def \hfdv { {\widehat f^2 } }

\def \fo {  f^1 }
\def \fdv { f^2  }

\begin{document}
\title{Hardy Spaces on Compact Riemann Surfaces with Boundary}
\author{
A. Zuevsky
 \\
School of Mathematics,  
Statistics and Applied Mathematics,  
\\
National University of Ireland, 
Galway, Ireland
}
\maketitle
\begin{abstract}
\noindent
We consider the holomorphic unramified  mapping of
 two arbitrary finite bordered Riemann surfaces.
Extending the map to the doubles $X_1$ and $X_2$
of Riemann surfaces we define the vector bundle on the second double
as a direct image of the vector bundle on first double. 
We choose line bundles of half-order differentials
 $\Delta_1$ and $\Delta_2$
so that the vector bundle  $V^{X_2}_{\chi_\ad}
 \otimes \Delta_2$ on $X_2$ would
be  the direct image of the vector bundle
 $V^{X_1}_{\chi_\ao} \otimes \Delta_1$.
We then  show that the Hardy spaces
$H_{ 2, J_1(p) } (S_1,V_{\chi_\ao} \otimes \Delta_1)$
and $H_{ 2,J_2(p) } (S_2,V_{\chi_\ad} \otimes \Delta_2)$
are isometrically isomorphic.
 Proving that we construct an explicit
isometric isomorphism and a matrix representation $\chi_\ad$
of the fundamental group $\piod(X_2, p_0)$ given
a matrix representation $\chi_\ao$ of the fundamental group
$\piod(X_1, p'_0)$.
On the basis of the results of \cite{vin} and Theorem \ref{theorem_1} 
proven in the present  work
we then conjecture that there exists a covariant functor
from the category  ${\cal RH}$ of finite bordered surfaces
with vector bundle and signature matrices
to the category of Kre\u{\i}n spaces and isomorphisms
which are ramified covering of Riemann surfaces.
\end{abstract}

\newpage
\section{Introduction}
It is well known how to study Hardy spaces defined on a finite bordered
Riemann surface \cite{abr}, \cite{rudin}, \cite{sarason}, \cite{voi1}.
For  domains  with more then one boundary component
it is natural to introduce, besides the usual positive definite inner
product on $H_2$, indefinite inner products. Those products may be
 introduced
by picking up different signature matrices
when integrating over different components of
the boundary of the Riemann surface.
In the paper \cite{vin}
a necessary and sufficient condition for such an indefinite
inner product to be non-degenerate was obtained.
 It was shown
that when this condition is satisfied one actually
gets a Kre\u{\i}n space. 
The result was obtained by using a covering map of the surface to the unit disk to 
construct an isomorphism to a Hardy-Kre\u{\i}n space over the unit disk. 
 Furthermore, each holomorphic mapping of
the finite bordered Riemann surface onto the unit disk
(which maps boundary to boundary)
determines an explicit isometric
isomorphism between this space and a usual vector-valued Hardy space
on the unit disk with an indefinite inner product defined by
an appropriate hermitian matrix.
The mapping to the unit disk in \cite{vin} serves as a tool 
to study of the Hardy-Kre\u{\i}n space over the finite bordered Riemann surface which in 
turn has motivation from the point of view of the study of commuting tuples of non-selfadjoint 
operators. 
As it is usual when studying Hardy spaces
on a multiply connected domain, the elements of the space are sections
of a vector bundle rather than functions.
The main point of the paper \cite{vin} was to construct
an appropriate extension of this bundle to the double of
the finite bordered Riemann surface and to use Cauchy kernels for certain
vector bundles on a compact Riemann surface.
Hardy spaces on a finite bordered Riemann surface, including
indefinite Hardy spaces, are important in the model theory
for commuting non-selfadjoint operators \cite{vino}.
 
Half-order differentials play a very important role in the vertex operator 
algebra approach to construction of partition and $n$-point functions 
for conformal field theories defined on Riemann surfaces \cite{mtz1}, \cite{mtz11}, \cite{mtz2}. 
In particular, the Szeg\"o kernel \cite{fay2} 
 turned out to be key object in construction of correlation 
 functions in free fermion conformal field theories/chiral algebras on 
a genus two Riemann surface sewed from two genus one Riemann surfaces \cite{y}. 

 In this work we replace a holomorphic mapping of a finite bordered
Riemann surface onto the unit disk by a holomorphic mapping
 of two arbitrary finite bordered
Riemann surfaces $S_1$ and $S_2$,
 which we assume however to be unramified.
In the spirit of \cite{vin} (see also \cite{bali}) one can introduce the extension of the
vector bundles on a finite bordered Riemann surfaces
to the  respective  doubles.
Extending the map $F$ to the doubles $X_1$ and $X_2$ of Riemann surfaces
$S_1$ and $S_2$ we define the vector bundle
$V^{X_2}_{\chi_\ad}$ on $X_2$
as a direct image of the vector bundle $V_{ \chi_\ao}^{X_1}$ over $X_1$.
We choose line bundles of half-order differentials
(i.e.,  square roots  of the canonical bundles $K_{X_i}$, $i=1,2$)
$\Delta_1$ and $\Delta_2$
so that the vector bundle  $V^{X_2}_{\chi_\ad}
 \otimes \Delta_2$ on $X_2$ would
be  the direct image of the vector bundle
 $V^{X_1}_{\chi_\ao} \otimes \Delta_1$.
We then  show that the Hardy spaces
$H_{ 2, J_1(p) } (S_1,V_{\chi_\ao} \otimes \Delta_1)$
and $H_{ 2,J_2(p) } (S_2,V_{\chi_\ad} \otimes \Delta_2)$ are isometrically
isomorphic. Proving that we construct a.) an explicit
isometric isomorphism;
b.) a matrix representation $\chi_\ad$
 of the fundamental group $\piod(X_2, p_0)$ given
a matrix representation $\chi_\ao$ of the fundamental group
$\piod(X_1, p'_0)$.
Using results of \cite{vin} and Theorem \ref{theorem_1} 
proven in the present  work we then conjecture that there exists a covariant functor
   from the category  ${\cal RH}$ of finite bordered surfaces
with vector bundle and signature matrices
to the category of Kre\u{\i}n spaces and isomorphisms
which are ramified covering of Riemann surfaces.

The isomorphism established  in this work
has also an operator theoretical interpretation,
namely, a (ramified) covering $F$ allows us to construct
a pair of commuting non-selfadjoint operators with the model
space on $S_2$ given  a pair of commuting non-selfadjoint operators
with the model space on $S_1$.
More generally, one might expect also possible connections with vessels construction and Bezoutians 
 \cite{kra}. 

We would like to thank V.Vinnikov for numerous discussions. 
\section{ Preliminaries}
As we mentioned in Introduction
 indefinite {\it Hardy spaces} \cite{hof}, \cite{dur}, \cite{garn}
on a finite
 bordered Riemann surface were considered in \cite{vin}.

\medskip
\noindent
{\bf Definition.} Let $J$ be an $m \times m$ unitary
self-adjoint matrix. Such a  matrix is
usually called  a  {\it signature matrix}.
(In fact one may take $J$ to be any non-singular self-adjoint matrix).
The Hardy space $H_2^m$
on the unit disk {\bf D}
 endowed with the indefinite
inner product
\begin{equation*}
[f,g]_J=\frac{1}{2\pi} \int_0^{2\pi} g(e^{it})^* J f(e^{it})dt,
\end{equation*}
is a {\it Kre\u{\i}n space} denoted $H^m_{2,J}$.

\medskip
\noindent
This space plays an important role in interpolation theory \cite{dym}
and in model theory  \cite{bra}.
  For the general theory of Kre\u{\i}n spaces
see \cite{bogn}, \cite{ioh}, \cite{aziz}.

Suppose now that  we have an open Riemann surface $S$
 such that
$S \cup \partial S$ is
 a finite bordered Riemann surface (i.e., a compact Riemann
 surface  with boundary), with the  boundary $\partial S$
consisting of $k \ge 1$ components ${\cal X}_0, ..., {\cal X}_{k-1}$.
We consider analytic sections of a rank $m$ flat
 unitary vector bundle $V_\chi$ on $S$ corresponding to
  a homomorphism $\chi$
from the fundamental group $\pi_1(S, p_0)$ into the group  $U(m)$
of $m \times m$ unitary matrices. An analytic
section $f$ of $V_\chi$ over  $S$ is an analytic ${\mathbb C}^m$-valued function
on the universal covering $\widetilde{S}$ of $S$ satisfying
\begin{equation*}
f(T\tilp)=\chi(T)f(\tilp), 
\end{equation*}
for all $\tilp \in {\widetilde S}$ and all deck transformations
$T$ of ${\widetilde S}$ over $S$, which we identify with elements
of the fundamental group $\piod(S,p)$; $f$ can be thought of as a
 multiplicative multivalued function on $S$.
We consider also multiplicative half-order differentials \cite{vin}, i.e.,
sections of a vector bundle of the form $V_\chi \otimes \Delta $, where
$V_\chi$ is a flat unitary vector bundle on
 $S$ as above and $\Delta$ is a square  root  of the
 canonical bundle on $S$: $\Delta \otimes \Delta \cong K_S$.

\medskip
\noindent
{\bf Definition.} The {\it Hardy space}
 $H_2(S, V_\chi \otimes \Dlt)$ on a Riemann surface $S$
 is the set of sections $\hf$
of a vector bundle $V_\chi \otimes \Delta$ analytic in $S$
satisfying
\begin{equation}
\label{ff}
\sup\limits_{1-\epsilon< r<1}
\sum_{i=0}^{k-1} \int_{\Xsi(r)} \hf(p)^*\hf(p) <  \infty, 
\end{equation}
for some $\epsilon >0$. In (\ref{ff})  $\Xsi(r)$ denotes smooth
 simple closed  curves
in  $S$ approximating $\Xsi$, $i=0, ... , k-1$ (\cite{vin}):
if $z_i$ is a boundary uniformizer near the boundary component $\Xsi$
then $\Xsi(r)$ is given by $|z_i(p)|=r$.

Note that since $\hf(p)$
 is a section of $V_\chi \otimes \Delta$, the expression
$\hf(p)^*\hf(p)$ is a section of $|K_S|$, where $|K_S|$ is
the line bundle with transition
functions the absolute values of the transition functions of $K_S$;
 sections  of $|K_S|$  can be represented locally as
 $\eta(t)|dt(p)|$ where $t(p) $ is a local parameter.
Therefore one can integrate  $\hf(p)^*\hf(p)$
over curves in $S$
and (\ref{ff}) makes sense.

The space $\htwoch$ is a Hilbert space with the inner product
\begin{equation*}
\langle \hf,\hg \rangle=\lim\limits_{r \lmpt 1} \sum_{i=0}^{k-1} \int_{\Xsi(r)} \hg(p)^*
\hf(p).
\end{equation*}
For a relation between $H_2(S, V_\chi \otimes \Delta )$ and Hardy
spaces of functions on $S$ with respect to a harmonic measure on
$\partial S$, see \cite{vin}.

\medskip
\noindent
{\bf Definition.} Denote by $H_{2, J(p)}(S, V_\chi \otimes \Dlt)$
the analogue of the
{\it Kre\u{\i}n} space
$H^m_{2,J}$
 for $S$
which is the Hardy space $H_2 (S, V_\chi \otimes \Dlt)$
 endowed with the indefinite inner product
\begin{equation}
\label{prod}
[\hf,\hg]_{J(p)}=\sum_{i=0}^{k-1} \int_{ {\cal X }_i } \hg(p)^* J(p)
 \hf(p), 
\end{equation}
where  $\hf(p)$ is the non-tangential boundary
values of $\hf(p)$ which exists again almost everywhere on $\partial S$
(see \cite{vin}) and
$J(\widetilde{p})$ is a locally constant matrix function
 on $\partial {\widetilde S}$
whose
values are $m \times m$ signature matrices, satisfying
\begin{equation}
\label{cond}
\chi(T)^*J(T\widetilde{p})\chi(T)=J(\widetilde{p}), 
\end{equation}
for
all $\widetilde{p} \in \partial \widetilde{S}$ and all $T \in \pi_1(S, p_0)$.
The expression $\hg(p)^* J(p) \hf(p)$ in (\ref{prod}) means
$\hg(\tilp)^* J(\tilp) \hf(\tilp)$ where $\tilp \in \partial {\widetilde S}$
is over $p \in \partial S$. It is a well defined section of $|K_S|$
because of the transformation property of $J(\tilp)$.
There exists a certain freedom in the choice of $J(\tilp)$ for the given
$V_\chi$. Indeed, choose points $p_i\in {\cal X}_i$, $i=0,...,k-1$.
Let $C_i$ be a path on $S$ linking $p_0$ to $p_i$. Set
$A_i=C_i^{-1}{\cal X}_i C_i \in \piod(S, p_0)$ (see Appendix).
Then the (homotopy class of)
$C_i$ determines a component ${\widetilde {\cal X} }_i$
of $\partial \widetilde{S}$ lying over ${\cal X}_i$,
and the constant value $J_i$ of $J(\tilp)$ on ${\widetilde {\cal X} }_i$
can be an arbitrary $m \times m$ signature matrix satisfying
\begin{equation*}
\chi(A_i)^* J_i \chi(A_i)=J_i. 
\end{equation*}
 Any other component of
 $\partial {\widetilde S}$ lying over ${\cal X}_i$ can be
obtained from
${\widetilde {\cal X} }_i$  by some deck transformation $R$. The
value of $J(\tilp)$ on this component is $\chi(R)^* J_i \chi(R)$.
For the case of the line bundles (i.e., $m=1$), the choice of $J(\tilp)$
amounts to an arbitrary choice of a sign $\pm 1$ for each ${\cal X}_i$.
We will often assume the choice of  components ${\widetilde {\cal X}}_i$
has been made and denote $H_{2,J(p)}(S, V_\chi \otimes \Delta)$
by $H_{2, J_0,..., J_{k-1}} (S, V_\chi \otimes \Delta)$.
The space $H_{2, J_0,..., J_{k-1}} (S, V_\chi \otimes \Delta)$ is
a natural example of an indefinite inner product space.
It is related to the model theory of pairs of
commuting non-selfadjoint operators and interpolation theory on
multiply connected domains.

In the paper \cite{vin} an appropriate extension of $V_\chi$ on $S$ to the
double $X$ of the Riemann surface $S$ was constructed.
Given a flat unitary vector bundle on a finite
bordered Riemann surface, together with a collection of signature
matrices, it can be uniquely extended to a flat
unitary vector bundle on the double satisfying certain symmetry properties.
Let us recall that construction.

 Due to the identification of the
boundaries the complex structures on two copies of $S$ constituting $X$
are mirror images of  each other, i.e., there exists an
{\it anti-holomorphic involution} $\tau: X \lngra X$
that maps ${\overline S}$ to ${\overline S'}$.
Thus $X$ is a compact real Riemann surface,
or equivalently
Riemann surface of a real algebraic curve.
The  genus $g$ of the double of $X$ of $S$ is
$g=2s+k-1$, where $s$ is the genus of $S$.
The set $X_f$ of fixed points of $\tau$ (real points of $X$) coincides
with the boundary $\partial S$ of $S$.
Furthermore $X$ is a real Riemann surface of dividing type:
the complement $X \backslash X_f$ consists of two connected
components $X_+ = S$ and $X_-=S'$ interchanged by $\tau$.
The converse is also true: any real Riemann surface of dividing
type is the double of a finite bordered Riemann surface.
 The anti-holomorphic involution  $\tau$
acts both on the fundamental group $\pi_1(X, p_0)$
 and on the universal covering
$\tilde{X}$ of $X$
(recall that the fundamental group $\piod(X, p_0)$ is isomorphic to
the group of deck transformations
${\it Deck}({\widetilde X}/X)$.
It also acts naturally on complex holomorphic vector
bundles on $X$: the transition
functions for the vector bundle $V^\tau$ are complex conjugates of the
transition functions for $V$ at the point conjugate under $\tau$.

Consider a vector bundle $H$ on $X$  of
  rank $m$  with ${\rm deg}~H=m(g-1)$
satisfying the condition $h^0(H)=0$.
 Such a vector bundle is  necessarily of the form
$H\cong V_\chi\otimes\Delta$ where $V_\chi$
is a rank $m$ flat vector bundle on
$X$ and $\Delta$ is a square root of $K_X$  \cite{balvin}.
  These vector
bundles are closely related to determinantal representations of algebraic
curves and play an important role in the theory of commuting non-selfadjoint
operators and related theory of 2D  systems
 \cite{vvv1}, \cite{vvv2}, \cite{vvv3}, \cite{balvin}, \cite{vvv4}

Let $H$ be such that there exists a  non-degenerate bilinear
pairing $H\times H^\tau\rightarrow K_X$ which is parahermitian.
The parahermitian property  means that
\begin{equation*}
(\hf,\hg^\tau)(p)=\overline{ (\hg,\hf^\tau)(p^\tau)}, 
\end{equation*}
for all local holomorphic
sections $\hf$ and $\hg$ of $H$ near $p$ and $p^\tau$ respectively.
We assume also
that the line bundle $\Delta$ has been chosen so that $\Delta\cong
\Delta^\tau$ and that  the transition functions of $\Delta$
are symmetric with respect to $\tau$ \cite{fay1};
Then we  obtain a
parahermitian non-degenerate bilinear pairing $V_\chi\otimes V_\chi^\tau
\rightarrow {\cal O}_X$, or more explicitly an everywhere nonsingular
holomorphic
$m\times m$ matrix--valued function $G$ on the universal covering $\tilde{X}$
with the property
\begin{equation}
\label{paraher}
G(\tilde{p}^\tau)^*=G(\tilde{p}), 
\end{equation}
satisfying the relation
\begin{equation}
\label{symrel}
\chi(T^\tau)^*G(T\tilde{p})\chi(T)=G(\tilde{p}), 
\end{equation}
where $T \in \piod(X, p_0)$.
The pairing $H\times H^\tau\rightarrow K_X$ is then given explicitly  by
\begin{equation*}
(\hf,\hg)(\tilde{p})=\hg(\tilde{p}^\tau)^*G(\tilde{p})\hf(\tilde{p}).
\end{equation*}

 Now let us   introduce the (in general) indefinite inner product
\begin{equation}
\label{bourbaki}
\left[ \hf, \hg \right]_{G(p)}  =
\int_{  X_f } \hg(\tilde{p}^\tau)^*G(\tilde{p})\hf(\tilde{p}), 
\end{equation}
where $\hf$ and $\hg$ are measurable sections of $H$ over $X_f$.
 Here and in similar expressions, the integral is computed on $X$ and the
integrand does not depend on the choice of $\tilde{p}\in\tilde{X}$ above
$p\in X$. Since in (\ref{bourbaki}) $\tilde{p}\in\tilde{X}$ lies over a point
of $X_f$ there exists $T_{\tilde{p}}\in\pi_1(X, p_0)$ such that
$\tilde{p}^\tau
=T_{\tilde{p}} \tilde{p}$. Therefore (\ref{bourbaki}) can be rewritten
as
\begin{equation}
\label{bourbaki1}
\left[\hf, \hg \right]_{ G(p) } =
\int_{  X_f }\hg(\tilde{p})^*J(\tilde{p})\hf(\tilde{p}), 
\end{equation}
where
\begin{equation*}
J(\tilde{p})=\chi(T_{\tilde{p}})^*G(\tilde{p}).
\end{equation*}
 Note that $J(\tilde{p})^*=J(\tilde{p})$ and
\begin{equation}
\label{symk}
\chi(R)^*J(R\tilde{p})\chi(R)=J(\tilde{p}), 
\end{equation}
 for all $\tilde{p}\in\tilde{X}$
over $X_f$ and all $R\in\pi_1(X, p_0)$. Thus the vector
bundle $H=V_\chi\otimes\Delta$ on $X$ defines an indefinite inner product
on the sections of its restriction to $X_f=\partial S$.

For $\tilde{p}$ lying over a point of $\Xsi$, we have
$T_{\tilde{p}}=R^\tau R^{-1}$ for $i=0$ and
$T_{\tilde{p}}=R^\tau B_iR^{-1}$ for $i=1,\dots ,k-1$
where $B_i=(C^\tau_i)^{-1} C_i$ are
 part of the
generators of the fundamental group
$\piod(X, p_0)$ of $X$ (see Appendix for the relation
between generators of the fundamental groups of a Riemann surface $S$ and
 the corresponding double), and
$R$ depends
only on the component of the inverse image of $\Xsi$ in $\tilde{X}$ that
$\tilde{p}$ belongs to. Restricting $\tilde{p}$ in (\ref{bourbaki1}) to belong
to a specific component we may write
\begin{equation*}
J(\tilde{p})=G(\tilde{p}), 
\end{equation*}
 for $i=0$ and
\begin{equation*}
J(\tilde{p})=\chi(B_i)^*G(\tilde{p}), 
\end{equation*}
 for $i=1,\dots ,k-1$ (see Appendix). (The 
specific component depends on the choice of the generators $B_i$, i.e.,
on the homotopy classes of the paths $C_i$.)

 It follows  from the conditions
${\rm deg}~H=m(g-1)$ and $h^0(H)=0$ that $H$ is a semi-stable vector
bundle.
By a theorem of Narasimhan and Seshadri \cite{sesh}
 $H$ is a direct sum of stable bundles if and only if the flat vector
bundle $V_\chi$
(in $H=V_\chi\otimes\Delta$) can be taken to be unitary flat. Since $G$ is an
isomorphism from $V_\chi$ to the dual of $V_\chi^\tau$ it follows in this
case that $G$ is constant and unitary. Since it is also selfadjoint, it
is a constant signature matrix. Thus  for the analytic sections
$\hg$ and $\hf$ of $V_\chi \otimes \Delta$
on $S$ that belong to $H_2(S, V_\chi \otimes \Delta)$
we can rewrite the inner product
 (\ref{bourbaki1})
as
\begin{equation}
\label{prod2}
\begin{array}{l}
\dsps
\bigskip
[\hf,\hg]_{ G(p) }
= \sum_{i=0}^{k-1}\int_{\Xsi} \hg(\tilp)^*J_i \hf(\tilp)
\\
\qquad \qquad
= [\hf,\hg]_{H_{2,J_0,...,J_{k-1}}(S, V_\chi \otimes \Delta) }, 
\end{array}
\end{equation}
where
\begin{equation*}
J_0=G, \qquad J_i=\chi(B_i)^*G, 
\end{equation*}
for $i=1,\dots ,k-1$ and $p$ is
restricted to belong to a specific component of the inverse image of $\Xsi$
in $\tilde{X}$ as explained above.
We then obtain for the vector bundle $H$ on $X$ the inner product
(\ref{prod}) on the Hardy space
$H_{2,J_0,\dots ,J_{k-1}}(S,V_\chi\otimes\Delta)$.
Conversely, every unitary flat vector bundle
on $S$ with signature matrices $J_0,\dots ,J_{k-1}$
 can be obtained from
a vector bundle $H$ on $X$ as above. 
Let $J_i$, $i=0,\dots ,k-1$ be
selfadjoint matrices and let $\chi:\pi_1(S, p_0)
\rightarrow {\rm GL}(m,{\mathbb C})$
be a homomorphism satisfying
\begin{equation*}
\chi(A_i)^*J_i\chi(A_i)=J_i.
\end{equation*}
Then by Proposition 2.1 from \cite{vin} there exists a
unique extension (still denoted by $\chi$) of $\chi$ to a homomorphism
from $\pi_1(X, p_0)$ into ${\rm GL}(m,{\mathbb C} )$ satisfying
\begin{eqnarray*}
\label{plkj}
\chi(T^\tau)^*G\chi(T)&=& G,\qquad T\in \pi_1(X),\\
\label{plkjh}
\chi(B_i)^*G&=&J_i
\end{eqnarray*}
where $G=J_0$ (see Appendix).

  If the original flat vector bundle
$V_\chi$ on $S$ is unitary flat and all the  matrices $J_i$
 are unitary then the extended vector bundle is also unitary flat
as it is follows form the proof of Proposition 2.1 of \cite{vin}.
 The extension need not satisfy $h^0(X,V_\chi\otimes\Delta)=0$; i.e.,  the
unitary case, this condition will be satisfied "generically" since
flat unitary vector bundles $V_\chi$ on $X$ with
$h^0(X,V_\chi\otimes\Delta)>0$
form a divisor in the moduli space of flat unitary vector bundles
(the generalized theta divisor \cite{drez}, \cite{fay2}).
It was proven in Proposition 2.2 from \cite{vin} that 
if the indefinite inner product space
$H_{2,J_0,\dots ,J_{k-1}}(S,V_\chi\otimes\Delta)$
is non-degenerate, then
\begin{equation}
\label{cocond}
h^0(X,V_\chi\otimes\Delta)=0.
\end{equation}
It follows that the  condition (\ref{cocond})
is satisfied automatically in the positive definite case
(i.e., when $J_i>0$ for $i=0,\dots ,k-1$);
for line bundles, this has been obtained in \cite{fay1}
(see also \cite{vvv5}).

Summing up, we see that the above extension procedure
establishes a one-to-one correspondence between
unitary flat vector bundles on $S$ together with a choice of
signature matrices satisfying (\ref{cond}) and
unitary flat vector bundle on $X$ satisfying the symmetry condition
(\ref{paraher}), (\ref{symrel}).
Given a unitary flat vector bundle on $S$, the
various choices of  extension
to  the double $X$ correspond to the various choices
of signature matrices.  We shall occasionally denote the corresponding
unitary flat vector bundles on $S$ and $X$  by $V_\chi^S$ and
$V_\chi^X$ respectively.

Under the condition $h^0(X, V_\chi \otimes \Delta)=0$  (i.e.,
that $V_\chi \otimes \Delta$ has no global holomorphic sections),
it turns out that $V_\chi \otimes \Delta$
 on $X$ admits a certain kernel function
(which is called the Cauchy kernel) which is an analogue of $\frac{I_m}{z-w}$
for the trivial bundle on the complex plane.
The Cauchy kernel
is the reproducing kernel for $H^2_{J(p)}(S, V_\chi \otimes \Delta)$.
In the case of line bundles the Cauchy kernel can be given explicitly
in terms of theta functions \cite{fay1}, \cite{ball}.
In \cite{vin} the Cauchy kernel was used to construct
for any given holomorphic mapping $z: S \longrightarrow {\bf D}$
an explicit
isometric isomorphism between  $\hojch$ and $H^M_{2,J}$
 for appropriate $M$ and $J$. In particular this implies that
$H_{2,J(p)}(S, V_\chi \otimes\Delta)$ is indeed non-degenerate
(under the condition $h^0(X, V_\chi \otimes \Delta)=0$) and actually a
Kre\u{\i}n space.
\section{ Statement of the Main Result}
 Suppose that we have two finite bordered Riemann surfaces $S_1$ and $S_2$.
Let $F: S_1 \longrightarrow S_2$ be an analytic mapping continuous up to the
boundary.
 Equivalently we may take $F$
to be a complex analytic mapping
$F: X_1 \longrightarrow X_2$ between the doubles of $S_1$ and $S_2$
equivariant with respect to the action of the anti-holomorphic involutions,
 i.e., such that the diagram
\begin{equation*}
\label{diatau}
\begin{array}{ccc}
\bigskip
X_1 & {\tau_1 \atop \longrightarrow }& X_1 \\
\bigskip
F \downarrow & & \downarrow F \\
\bigskip
X_2 &{ \tau_2 \atop  \longrightarrow} & X_2 
\end{array}
\end{equation*}
is commutative. Notice that $F: S_1 \lngra S_2$ is unramified if and
only if $F: X_1 \lngra X_2$ is unramified.

We identify as usual a complex holomorphic vector bundle
on a complex manifold with a locally free sheaf of its
analytic sections. It is easily seen that if $V^X$
is a complex holomorphic vector bundle of rank $m$
on a complex manifold $X$ and $F$ is a $n$-sheeted unramified covering,
then  the direct image $V^Y=F_* V^X$
is a complex holomorphic vector bundle of rank $nm$ on $Y$
and the fiber of $V^Y$ at a given point of $Y$
 is the direct sum of the fibers of $V^X$ at the preimages
of this point on $X$.

The main statement of this work is the following
\begin{theorem}
\label{theorem_1}
Let $F: S_1 \lngra S_2$ be a map of  finite bordered Riemann
surfaces which is a finite $n$-sheeted unramified covering $(F; S_1, S_2)$,
and
let $J_1(\tilp)$  be signature matrices
for a unitary flat vector bundle $V_{\chi_\ao}$ on $S_1$ of rank $m$.
Consider the corresponding  extension of
$V_{\chi_\ao}$ to the
double $X_1$ of $S_1$ satisfying the symmetry condition
\begin{equation*}
\chi_\ao (R^\tau)^*\; G_1(R\tilp)\; \chi_\ao (R)= G_1(\tilp), 
\end{equation*}
for all $R\in \piod(X_1, p_0')$ and all $\tilp \in {\widetilde X}$.
Choose the bundles $\Delta_1$ and $\Delta_2$
of half-order differentials on $X_1$ and $X_2$ respectively,
such that

\medskip
a.) the bundles $\Delta_i$, $i=1,2$ are  invariant with respect to the
corresponding anti-holomorphic involutions, i.e.,
$\Delta_i^{\tau_i}=\Delta_i$ and the transition functions of
$\Delta_1$ and $\Delta_2$ are symmetric with respect to $\tau_1$
and $\tau_2$;

\medskip
b.) the pull--back of  $\Dlt_2$ is  equal to  $\Dlt_1$, i.e.,
  $\Dlt_1=F^*\Dlt_2$.

\medskip
Then

\medskip
1.) the direct image  $V^{X_2}_{\chi_\ad}=F_* V_{\chi_\ao}^{X_1}$
is a unitary flat holomorpnic vector bundle of rank $nm$
satisfying the
symmetry condition
\begin{equation*}
\chi_\ad (T^\tau)^*\; G_2(T\tilp)\; \chi_\ad (T)=G_2(\tilp), 
\end{equation*}
for all $T\in \piod(X_2, p_0)$ and all $\tilp \in {\widetilde X}$,
and an appropriate matrix function $G_2(\tilp)$ and representation
$\chi_\ad$ of $\piod(X_2, p_0)$; furthermore $F_*(V_{\chi_\ao} \otimes
\Delta_1)= V_{\chi_\ad} \otimes \Delta_2$;

\medskip
2.) there exists a canonical isometric
isomorphism
\begin{equation*}
\label{isoso}
\phi_F: H_{2,J_1(p)}(S_1, \vchio \otimes \Delta_1)\;
\widetilde{ \longrightarrow }\;
  H_{2,J_2(p)}(S_2, \vchidv \otimes  \Delta_2), 
\end{equation*}
between Hardy spaces on $S_1$ and $S_2$.
\end{theorem}

\medskip
Now some remarks are in order.
By definition the anti-holomorphic involutions $\tau_1$ and
$\tau_2$ are related by
\begin{equation*}
F \circ \tau_1 = \tau_2 \circ F, 
\end{equation*}
and therefore if we have a
line bundle $L_2$ on $X_2$ then its  pull-back satisfies
\begin{equation*}
\left( F^* L_2^{\tau_2}\right) = \left( F^* L_2\right)^{\tau_1}.
\end{equation*}
 We fix $\Delta_2$ such that $\Delta_2^\tau=\Delta_2$ and
$\Delta_1=F^*\Delta_2$. Then it follows that $\Delta_1^\tau=\Delta_1$.
We choose $\Delta_2$ such that its
transition functions  are symmetrical. Then since
$F$ is equivariant with respect to the anti-holomorphic involution then
 transition functions of $\Delta_1$ are also symmetrical.

The isomorphism of the spaces
 $H_{2, J_1(p)}(S_1, V_{\chi_\ao}  \otimes \Dlt_1)$
and $H_{2, J_2(p)}(S_2, \vchidv \otimes \Dlt_2)$ implies
that they are degenerate or non-degenerate simultaneously,
 i.e.,
 $h^0(X_1, \vchio \otimes \Delta_1)=h^0(X_2, V_{\chi_\ad}
 \otimes \Delta_2)= 0$ which is obvious from the definition
of the direct image vector bundle.

We assume that
the map $F: S_1 \lmpt S_2$ is a $n$-sheeted {\it unramified} covering
$(F; S_1, S_2)$ of the Riemann surface $S_2$ by $S_1$.
On the other hand, a result of \cite{vin} mentioned
in Introduction is a construction of an isometric
isomorphism between Hardy spaces when $S_2={\bf D}$ but $F$
is (usually) ramified  ( assuming
$H_{2, J_1(p) (S_1, V_{\chi_\ao} \otimes \Delta_1 ) }$
is not degenerate, i.e., $h^0(X_1,V_{\chi_\ao} \otimes \Delta_1)=0$).
 The next natural step would be to consider the case when $S_2$
is an arbitrary finite bordered Riemann surface and $F$ is a
ramified covering. That will be a point of some further
publication.

We have introduced the vector bundle $V^{X_2}_{\chi_\ad}$  on the double $X_2$
as the direct image of the vector bundle $V_{\chi_\ao}^{X_1}$ on $X_1$
defining the vector bundle $V^{S_2}_{\chi_2}$ on $S_2$ and the signature
matrices $J_2(\tilp)$.
On the other hand one
can define  the vector bundle $V^{S_2}_{\chi_\ad}$ to be the direct
image $F_* V_{\chi_\ao}^{S_1}$ of the vector bundle $V_{\chi_\ao}^{S_1}$
 with signature matrices defined naturally in terms
of $J_1(\tilp)$ (as direct sums)
 Though the main claim of Theorem \ref{theorem_1} 
is formulated  for finite bordered
Riemann surfaces it seems to us that the
consideration of the structures involved in its proof is more
natural (in the sense of the theory of compact Riemann surfaces)
on the doubles. Furthermore, this approach allows us to
construct a matrix representation $\chi_\ad$ of the
fundamental group $\piod(X_2, p_0)$ given a representation of
$\piod(X_1, p_0')$, and the matrix function $G_2(\tilp)$.
We will prove that signature matrices $J_2(\tilp)$ calculated
with the help of the representation $\chi_\ad$
do coincide with the signature matrices constructed
directly from the signature matrices $J_1(\tilp)$.
This shows the equivalence of those two approaches.
From the use of Cauchy kernels in \cite{vin} it seems
however that in the ramified case the approach via the
doubles is the only one possible.

Speaking in more abstract terms
 we deal in Theorem \ref{theorem_1} with a category
which we will denote by ${\cal RH}$.
Objects of ${\cal RH}$  are
finite bordered Riemann surfaces $S$
 together with a
unitary flat vector bundle $V_\chi$ and
signature matrices $J(\tilp)$
(or equivalently,
compact real Riemann surfaces $X$ of dividing type with a
vector bundle $V_\chi^X \otimes \Delta$ on $X$ and
a matrix function $G(\tilp)$ satisfying (\ref{paraher})
 and (\ref{symrel}))
such that the space $H_{2,J(p)}(S, V_\chi \otimes \Delta)$
is non degenerate,
i.e., $h^0(X, V_\chi \otimes \Delta)=0$.
A morphism
between the objects
$(X_1, V_{\chi_\ao}^{X_1} \otimes \Delta_1, G_1)$
and  $(X_2, V_{\chi_\ad}^{X_2} \otimes \Delta_2, G_2)$ of ${\cal RH}$
is an analytic map $F: X_1 \lngra X_2$
 of Riemann surfaces which is equivariant with respect
to  anti-holomorphic involutions $\tau_1$ and $\tau_2$, such that
$V^{X_2}_{\chi_\ad}
 \otimes \Delta_2 = F_*(V_{\chi_\ao}^{X_1} \otimes \Delta_1)$
  (and $G_2(\tilp)$ is correspondingly induced by $G_1(\tilp)$).

\medskip
{\it We conjecture that there exists a covariant functor
from the above mentioned category ${\cal RH}$ to
 the category of Kre\u{\i}n spaces and isomorphisms,
 associating to $(S, V_\chi \otimes \Delta, J(\tilp))$
the Hardy space $H_{2,J(p)}(S, V_\chi \otimes \Delta).$}

\medskip
\noindent
Theorem \ref{theorem_1} proves the conjecture for a subcategory
of ${\cal RH}$ whose morphisms are unramified coverings.
 The isometric isomorphism established in \cite{vin} proves
another special case of the conjecture namely for a
subcategory whose morphisms are restricted to have  the
unit disk ${\bf D}$ as a range.
Somewhat related  considerations of categories
of functional spaces on
Riemann surfaces are contained  in \cite{all}.
\section{ Sections of the vector bundle $F_*(V_{\chi_\ao}\otimes \Delta_1)$ }
In this section we give an explicit construction of a
holomorphic section $\hfdv$ of the bundle
$F_*(V_{\chi_\ao}\otimes \Delta_1)$ on $X_2$ in terms of a
holomorphic section $\hfo$ of $\bundo$  on $X_1$.
When $(F;S_1, S_2)$ is an unramified covering the doubles
$X_1$ and $X_2$ possess the common universal covering ${\widetilde X}$,
i.e., one has  a diagram
\begin{equation*}
\label{kriv}
\begin{array}{ccc}
\dsps
 & {\widetilde  X}  &
\\
\dsps
{\pi_1 \atop {} } \swarrow  &   & \dsps \searrow  {\pi_2 \atop {} }
\\
\dsps
X_1 &  { F \atop \lngra  }  & X_2 
\end{array}
\end{equation*}
where $\piod$ and $\pi_2$ are the covering maps  from ${\widetilde X}$
to $X_1$ and $X_2$ respectively.
Let $U' \subset X_1$ be an open set in $X_1$.
Suppose $\fo$ is an analytic  section of the holomorphic vector
bundle $\vchio$ over $U'$, i.e.,
an analytic ${\mathbb C}^m$-valued function on $\piod^{-1}(U')
\subset {\widetilde X}$  satisfying the relation
\begin{equation}
\label{svo}
\fo(T\tilp)=\chi_\ao (T)\fo(\tilp).
\end{equation}
Similarly, a section $\hfo$
of the vector bundle $V_{\chi_\ao} \otimes \Delta$ over $U'$
satisfies
\begin{equation*}
\frac
{\hfo(T\tilp)}
{ \sqrt{dt_1(T \tilp)} }
=\chi_\ao (T)
\frac
{\hfo(\tilp)}
{ \sqrt{dt_1(\tilp)} }, 
\end{equation*}
for all $\tilp \in \piod^{-1}(U')$, $T \in \piod(X_1, p'_0)$, where
$t_1$ is a local parameter on $X_1$ lifted to ${\widetilde X}$.
 The fundamental group
$\piod(X_1, p'_0)$, $p' \in X_1$, is a subgroup of $\piod(X_2, p_0)$
of index $n$ (here $p'_0$ is a preimage of $p_0 \in X_2$).
Enumerate fixed representatives $g_i$, $i=1,...,n$ of the
left cosets $\left\{  \piod(X_1, p_0') g_i \right\}$ of the group
$\piod(X_2, p_0)$ with respect to its subgroup $\piod(X_1, p')$.
We define a sheaf on $X_2$ whose sections over an open set $U \subset X_2$
 are analytic ${\mathbb C}^{mn}$-valued functions on $\piod^{-1}(U)$ of
the vector form
\begin{equation}
\label{sect0}
\fdv(\tilp)=
\left[
\fo (g_i \tilp)
\right], 
\end{equation}
$i=1, \ldots, n$, 
where $\fo$ is a section of the bundle $V_{\chi_\ao}$
 over  $F^{-1}(U)$, i.e., $\fo(\tilp)$ is an analytic
${\mathbb C}^m$-valued function on $\piod^{-1}(F^{-1}(U))=\pi^{-1}_2(U)$
satisfying (\ref{svo}).
It easy to see from the definition
that this sheaf on $X_2$
 is isomorphic to the direct image sheaf of a sheaf on $X_1$
 of analytic sections of $V_{\chi_\ao}$, i.e., (\ref{sect0})
defines the sheaf of analytic sections of $F_*V_{\chi_\ao}$. 

Now let $p \in X_2$ and $p'_1, ..., p'_n \in X_1$ be preimages
of $p$. Let $t_2$ and  $t_{1,i}$ be
local parameters near $p$ and  $p'_i$, $i=1,...,n$ lifted to
the common universal covering ${\widetilde X}$.
Denote by $\varphi_i$ the composition $t_2^{-1} \circ F \circ t_{1,i}$.
 Then a section $\hfdv(\tilp)$
of a vector bundle  $F_*(V^{X_2}_{\chi_\ao} \otimes \Delta_1)$ is
given by
\begin{equation}
\label{sect1}
\hfdv(\tilp)=
\left[
\dsps
\frac{\hfo (g_i\tilp)}
 {\sqrt{\varphi'_i (g_i\tilp) }  \sqrt{ dt_{1,i} (g_i \tilp) } }
\right]\; 
\sqrt{dt_2(\tilp)}, 
\end{equation} 
where $\hfo(\tilp)$ is a section of the vector bundle
$V_{\chi_\ao} \otimes \Delta_1$
 and $t_{1,i}$, $i=1,...,n$ are local
parameters in the vicinity of $g_i \tilp$.
Since we have chosen
 the bundles $\Delta_1$ and $\Delta_2$  of  half-order
differentials in (\ref{sect1}),
 the ambiguity in the sign of the square roots of
$\varphi'(g_i \tilp)$ in (\ref{sect1}) is global and
since we have assumed that $\Delta_1=F^*\Delta_2$,
the expression (\ref{sect1}) does not depend on the choice of
local parameters.
\section{Representation $\chi_\ad$ of $\piod(X_2, p_0)$}
In this section we give an explicit formula for
a unitary representation $\chi_\ad$ of $\piod(X_2, p_0)$
such that $V_{\chi_\ad}=F_*V_{\chi_\ao}$.
It follows from the previous section that we  have
to define $\chi_2$ so that
\begin{equation}
\label{gado}
f^2(T \tilp)= \chi_\ad (T) f^2(\tilp), 
\end{equation}
for every $f^2$ given by (\ref{sect0}).
Let $g \in \piod(X_2, p_0)$. Fix a preimage $p_0' \in X_1$ of $p_0$.
The element $g$ belongs to a coset of the fundamental
group $\piod(X_2, p_0)$ with respect to its subgroup
 $\piod(X_1, p_0')$.
Then there exist elements $h \in \piod(X_1, p'_0)$ and
$g_{\sigma_g (i)} \in \piod(X_2, p_0)$
  such that
\begin{equation*}
g_i g =h g_{\sigma_g} (i),
\end{equation*}
i.e., $g$ defines a permutation $\sigma_g$ of the preimages
 of $p_0$. We take this as a {\it definition}  of $\sigma_g$.
We define the matrix representation $\chi_\ad$ as follows:
\begin{equation}
\label{repr2}
\left[
\chi_\ad (g)
\right]_{kj}
=
\chi_\ao (g_k g g^{-1}_{ \sigma_{g} (k)} )\; \delta_{ \sigma_{g} (k) j }.
\end{equation}
It is immediate that (\ref{gado}) is verified.
 Taking into account the unitarity of $\chi_\ao (g)$,
it can be seen from (\ref{repr2}) that  the matrices
defining  the representation of $\piod(X_2, p_0)$
are unitary, i.e.,
\begin{equation*}
\left[ \chi_\ad (g_i) \cdot  \chi_\ad (g_i)^* \right]_{kj}
= \delta_{kj}.
\end{equation*}
Now we check that (\ref{repr2}) provides a representation $\piod(X_2, p_0)$,
 i.e.,
\begin{equation*}
\chi_\ad (g \widetilde{g} )=\chi_\ad (g )\chi_\ad ( \widetilde{g} ), 
\end{equation*}
for all $\;g$, ${\widetilde g} \in \piod(X_2, p_0)$. Proving this
we used the fact that
 $\chi_\ao$ is homomorphism and
$\sigma_g$ is an anti-homomorphism, i.e.,
\begin{equation*}
\sigma_{g'}(\sigma_{g''}(k) )= \sigma_{g'' g' }(k),
\end{equation*}
which can be easily verified.
In general, the matrix $\chi_\ad$ is given by the formula
\begin{equation*}
\label{repr1}
\left[
\chi_2(g)\right]_{kj}= \chi_\ao (g_k g g_j^{-1})
\; \delta_{k+i-j-1,0}, 
\end{equation*}
where $g$ belongs to $i$-th coset.
\section{ Construction of pairing and  inner product }
 Suppose that
$H^1=V_{\chi_\ao}\otimes \Delta_1$ is such that there exists a
non-degenerate bilinear pairing $H^1 \times (H^1)^\tau \lngra K_{X_1}$
which is parahermitian, i.e.,
\begin{equation*}
\left(\hfo, \hgo^\tau \right)(p)=
\overline{\left(\hgo , (\hfo)^\tau  \right) (p^\tau) }.
\end{equation*}
We assume that
the line bundle $\Delta_1$ is such that
$\Delta^{\tau_1}_1 \cong \Delta_1$ and the transition functions
of $\Delta$ are symmetric with respect to $\tau_1$.
Then we have a parahermitian non-degenerate bilinear
pairing $V_{\chi_\ao} \otimes V_{\chi_\ao}^\tau \lngra {\cal O}_{X_1}$
 and the matrix function $G_1(\tilp)$ satisfying (\ref{paraher})
and (\ref{symrel}).
One can define a bilinear non-degenerate pairing
$H^2 \times (H^2)^\tau \lngra K_{X_2}$
where $H^2=F_*H^1=V_{\chi_\ad} \otimes \Delta_2$,
 introducing an everywhere
nonsingular holomorphic $mn \times mn$ matrix-valued function
$G_2$ on the universal covering $\widetilde{X}$ of $X_1$ and $X_2$.
The matrix $G_2(\tilp)$ should have the property
\begin{equation}
\label{prop3}
G_2(\tilp^\tau)^*= G_2(\tilp), 
\end{equation}
and satisfy the symmetry condition, $T \in \piod(X_2, p_0)$
\begin{equation}
\label{gedv}
\chi_2(T^\tau)^* G_2(T\tilp) \chi_2(T)= G_2(\tilp).
\end{equation}
Then the  pairing is given by
\begin{equation*}
\label{prodva}
\left(
\hfdv, \hgdv \right)(\tilp)=
\hgdv(\tilp^\tau)^* G_2(\tilp) \hfdv(\tilp).
\end{equation*}
Taking into account the explicit form (\ref{repr2}) of
$\chi_2(g)$ one can check that the following expression for $G_2(\tilp)$
does satisfies (\ref{gedv})
\begin{equation}
\label{gedvop}
\left [ G_2(\tilp) \right]_{kj}=
G_1(g_k^\tau \tilp)\; \chi_\ao (h_k)\;  \delta_{\nu(k),j}, 
\end{equation}
where $\nu(k)$ is defined as follows.
Consider the action of $\tau$ on an element $g\in \piod(X_2, p_0)$.
By definition we have $g^\tau=\tau g \tau^{-1}$.
For any $g_k$ that belongs to $k$-th coset of $\piod(X_2, p)$
with respect to $\piod(X_1, p'_0)$ there exist $h_k \in \piod(X_1, p'_0)$
and $g_{\nu(k)} \in \piod(X_2, p)$ such that
\begin{equation}
\label{defnu}
g_k^\tau = h_k g_{\nu(k)}.
\end{equation}
We define $\nu(k)$ by  (\ref{defnu}).
 One can check directly that (\ref{gedvop}) does satisfy conditions 
(\ref{prop3}) and (\ref{gedv}).

We saw in Introduction how to define an indefinite inner product
(\ref{prod}),(\ref{prod2})  on the Hardy space $\hjb$ using signature matrices
$J_0,...,J_{k-1}$.
Suppose that we have such an inner product on $\hjo$
\begin{equation*}
[\hfo,\hgo]_{\hjo}
= \sum_{i=0}^{k-1}\int_{ {\cal X}_{1,i} } \hgo(\tilp)^*J_{1,i} \hfo(\tilp).
\end{equation*}
Then we define an indefinite inner product on $\hjdv$
\begin{equation}
\label{prod22}
[\hfdv,\hgdv]_{\hjdv}
=\int_{  X_{2,f} } \hgdv(\tilp^{\tau_2})^* G_2(\tilp) \hfdv(\tilp).
\end{equation}
By the same reasons as in Introduction we can rewrite (\ref{prod22})
as
\begin{equation}
\label{prod23}
[\hfdv,\hgdv]_{\hjdv}
=\int_{ X_{2,f} } \hgdv(\tilp)^* J_2(\tilp) \hfdv(\tilp),
\end{equation}
where
\begin{equation*}
J_2(\tilp)=\chi_2 (T_{\tilp})^*G_2(\tilp), 
\end{equation*}
and introduce the matrices
\begin{equation}
\label{jdvaper}
J_{2,0}=G_2, \qquad J_{2,i}=\chi_2(B_{2,1})^*G_2, 
\end{equation}
where $B_{2,1} \in \piod(X_2, p)$ (see Appendix).
As in \cite{vin} the extension of the bundle $V^{S_2}_{\chi_2}$ on
the Riemann surface $S_2$ to the double $X_2$ depends on the
choice of the signature matrices $J_{2,0},...,J_{2,k-1}$
given by (\ref{jdvaper})
and  which satisfies the
symmetry condition (\ref{symk}).
On the other hand, one can define the signature matrix $J_2(\tilp)$
using the signature matrix $J_1(\tilp)$. One should have
\begin{equation}
\label{jdvapro1}
\chi_2(T)^* J_2(T \tilp) \chi_2(T) = J_2(\tilp), 
\end{equation}
for all $T \in \piod(X_2, p)$ and
\begin{equation}
\label{jdvapro2}
J_2(\tilp)^*=J_2(\tilp), 
\end{equation}
for all $\tilp \in {\widetilde X}$  over $p \in X_{2,f}$.
The matrix $J_2(\tilp)$ in the form
\begin{equation}
\label{jdva}
\left[ J_2(\tilp) \right]_{kj}=J_1(g_k \tilp) \; \delta_{kj}, 
\end{equation}
satisfies (\ref{jdvapro1}) and (\ref{jdvapro2}).
Then we check the commutativity of the diagram
\begin{equation*}
\label{diagr}
\begin{array}{ccc}
\bigskip
V_{\chio}^{S_1} & {ext,\; J_1 \atop \longrightarrow }& V_{\chio}^{X_1} \\
\bigskip
F^S_* \downarrow & & \downarrow F^X_* \\
\bigskip
V_{\chidv}^{S_2} &{ext, \; J_2 \atop  \longrightarrow} & V_{\chidv}^{X_2} 
\end{array}
\end{equation*}
where $ \;ext, J_i\; $ means the extension of the vector bundle
$V^{S_i}_{\chi_{ { {} \atop i  } } }$ on $S_i$ to the double
$X_i$.
I.e., we will show
  that the matrix $J_2(\tilp)$ defined by (\ref{jdvaper}) coincides
with (\ref{jdva}). It easy to check that
\begin{equation*}
\label{tilppro}
T_{R \tilp} R = R^\tau T_\tilp, 
\end{equation*}
for all $R\in \piod(X_2,p)$, $p \in X_{2,f}$ and
$T_\tilp \in \piod(X_2,p)$ such that $\tilp^\tau=T_\tilp \; \tilp$
where $\tilp$ lies over $p$. Using that we arrive at
\begin{equation*}
\dsps
\left[J_2(\tilp)\right]_{kj}= \chi_2(\ttilp)^*G_2(\tilp)
= J_1(g_k \tilp )  \delta_{kj}.
\end{equation*}
\section{ Proof of the isometricity}
We have constructed explicitly  a section $\hfdv$ (\ref{sect1})
of the bundle $V_{\chi_\ad}\otimes \Delta_2$ in terms of
a section $\hfo$ of the bundle $V_{\chi_\ao}\otimes \Delta_1$.
Now we will prove that the map $\hfo \longmapsto \hfdv$ is
 an isometric isomorphism of the space $\hso$ on the space
$\hsdv$. First let us show that $\hfo \in \hso$ if and only if
$\hfdv \in \hsdv$.

 Suppose $\hfo$ is a section of the bundle $\bundo$ and
$\hfo \in \hso$.
That means that $\hfo \in H_2(S_1, \vchio \otimes \Dlt_1)$,
i.e., $\hfo$ is an analytic in $X_1$ and
\begin{equation*}
\label{ffo}
\sup\limits_{1-\epsilon < r<1}
\sum_{i=0}^{k-1} \int_{{\cal X}_{1,i}(r)} \hfo(p)^*\hfo(p) <  \infty, 
\end{equation*}
for some $\epsilon$. Here ${\cal X}_{1,i}(r)$ are smooth simple curves
in $X_1$ approximating the $i$-th boundary of the $X_1$.
The space $\hso$ is the space
$H_2(S_1, \vchio \otimes \Dlt_1)$  endowed with the indefinite
inner product (\ref{prod})
\begin{equation*}
[\hfo,\hgo]_{\hsjo}
= \sum_{i=0}^{k-1}\int_{ {\widetilde  {\cal X} }_{1,i} }
\hgo(\tilp)^*J_{1,i} \hfo(\tilp).
\end{equation*}
Let ${\cal X}_{2,i}$ be a boundary component of $X_2$ and
 ${\cal X}_{1,i_j}$, $j=1,...,n_i$ be corresponding preimages on
$X_1$. The boundary uniformizer $z_1$ near the boundary
component is such that $z_1 p'_0= z_2 \circ F p'_0$. Then
the approximating curves ${\cal X}_{2,i}(r)$ are mapped to
the approximating curves ${\cal X}_{1,i_j}(r)$, $j=1,...,n_i$.
Due to the  construction given by the formula (\ref{sect1}) we see
that $\hfdv$ is an analytic and
\begin{equation}
\label{ffdv}
\begin{array}{l}
\bigskip
\dsps
\sum\limits_{i=0}^{k_2-1} \int_{{\cal X}_{2,i}(r)} \hfdv(p)^*\hfdv(p)
=
\sum\limits_{i=0}^{k_2-1}
\sum\limits_{j=1}^{n_i}
 \int_{{\widetilde {\cal X} }_{1,i_j }}
 \hfo(\tilp)^* \hfo (\tilp)
\\
\bigskip
\quad
\dsps
=
\sum\limits_{i=0}^{k_1-1}
\int_{{\cal X}_{1,i }}
\hfo(p)^* \hfo (p)
<  \infty.
\end{array}
\end{equation}
The summation in (\ref{ffdv}) with
upper limits $n_i$ is performed over the components ${\cal X}_{1,i_j}$,
$j=1,..., n_i$
that are preimages of ${\cal X}_{2,i}$.
Thus we infer  that $\hfdv$ belongs  to the
space $H_2(S_2, \vchidv\otimes \Dlt_2)$.
In the previous section we have introduced an indefinite inner product
in the space $H_2(S_2, \vchidv\otimes \Dlt_2)$. Thus we see that
a section $\hfdv$ of $\bundv$ constructed by the formula (\ref{sect1})
belongs to the space $\hsdv$.

  Finally, it remains to show that the inner product (\ref{prod23})
is isometric, i.e., that
\begin{eqnarray}
\label{raven}
& & 
\left[ \hfdv (\tilp), \hhdv (\tilp)  \right]_{\hsjdv}
\\
\notag
& & \qquad =
\left[ \hfo (\tilp), \hho (\tilp) \right]_{\hsjo}, 
\end{eqnarray}
where $\hfo$,$\hho$ and $\hfdv$,$\hhdv$ are sections of the vector bundles
$\bundo$ and $\bundv$ respectively.
Indeed,
consider the inner product of  two  sections of the
bundle $\bundv$
\begin{equation*}
\begin{array}{l}
\bigskip
{\displaystyle
\left[ \hfdv (\tilp), \hhdv (\tilp)  \right]_{\hsjdv}=
\sum\limits_{l=0}^{k_2-1}
\int_{ {\widetilde  {\cal X} }_{2,l}  }
 \hfdv ( \tilp )^* J_{2,l} \hhdv (\tilp)
}
\\
\bigskip
\qquad
{\displaystyle
= \int_{ {\widetilde {\cal X} }_{2,f} }
  \hfdv ( \tilp^\tau )^* G_2(\tilp)  \hhdv (\tilp)
}
\\
\bigskip
\qquad
{\displaystyle
=\int_{ {\widetilde {\cal X} }_{2,f} }
 \sum\limits_{i,j =1}^{n}
\hfo ( (g_i^\tau \tilp)^\tau )^* G_1 (g_i^\tau \tilp)
\hho (g_i^\tau \tilp)
\frac{dt_2 ( \tilp) }
{\varphi'_i (g_i\tilp ) dt_{1,i} (g_i\tilp ) }.
}
\end{array}
\end{equation*}
By the same argument that were used in the formulae (\ref{ffdv})
the last integral is equal to
\begin{equation*}
\begin{array}{l}
\bigskip
\dsps
\int_{ {\widetilde {\cal X}}_{1,f} }
\hfo (\tilp^\tau )^*  G_1(\tilp) \hho (\tilp)
= \sum\limits_{l=0}^{k_1-1}
\int_{{\widetilde   {\cal X}}_{1,l}  }
\hfo (\tilp)^* J_{1,l} \hho (\tilp)
\\
\bigskip
\qquad
\dsps
=\left[ \hfo(\tilp), \hho(\tilp)  \right]_{\hsjo}, 
\end{array}
\end{equation*}
where we use the invariance of sections of $\Delta_1$ with
respect to  deck transformations and the symmetry of
the their transition functions.
 Hence  we see that (\ref{raven}) holds.
That completes the proof of the isometricity.
\section{Appendix: Fundamental groups of $S$ and double $X$}
 Let us describe explicitly \cite{vin} the action of $\tau$
on the generators of $\pi_1(X, p_0)$.
Choose points $p_i\in \Xsi, i=0,\dots ,k-1$,
and let $C_i$ be a path on $S$ linking $p_0$ to $p_i$. Then
$\pi_1(S, p_0)$ is generated by
\begin{equation}
\label{gener1}
A_0,A_1,\dots, A_{k-1},A_1^\prime,B_1^\prime,\dots,
A_s^\prime,B_s^\prime,
\end{equation}
where $A_0={\cal X}_0$, $A_j=C_j^{-1}{\cal X}_jC_j$ for $j=1,\dots ,k-1$, and
$A_i^\prime$, $B_i^\prime$, $i=1,\dots,s$,
represent a canonical homology basis on
$S$ with the  intersection matrix
$\scriptsize{ \left(\begin{array}{cc}
0&I\\-I&0\end{array}\right)}$.
The generators of $\pi_1(S, p_0)$  satisfy a single relation
\begin{equation*}
\prod_{i=1}^sA_i^\prime B_i^\prime A_i^{\prime -1}
B_i^{\prime -1}\prod_{k-1}^0A_i=1.
\end{equation*}
Now consider the fundamental group
$\pi_1(X, p_0)$. It is generated by
\begin{equation*}
A_1,B_1,\dots ,A_{k-1},B_{k-1},A_1^\prime,B_1^\prime,\dots,A_s^\prime,
B_s^\prime,A_1^{\prime \prime},B_1^{\prime \prime},\dots,
A_s^{\prime \prime},B_s^{\prime \prime}.
\end{equation*}
The generators  $A_j$, $A_i^\prime$, $B_i^\prime$ are the same as in
(\ref{gener1})
\begin{equation*}
B_j=(C_j^\tau)^{-1}C_j, 
\end{equation*}
 for $j=1,\dots, k-1$
and
\begin{eqnarray*}
\label{toto1}
A_i^{\prime \prime}&=&B_i^{\prime \tau}, 
\\
\label{toto2}
B_i^{\prime \prime}&=&A_i^{\prime \tau}.
\end{eqnarray*}
The generators of $\pi_1(X, p_0)$
 satisfy a single relation  by  \cite{natan}
\begin{equation*}
\prod_{i=s}^{1}A_i^{\prime \prime}B_i^{\prime \prime}
A_i^{\prime \prime -1}B_i^{\prime \prime -1}\prod_{i=1}^s
A_i^\prime B_i^\prime A_i^{\prime -1}B_i^{\prime -1}
\prod_{j=k-1}^1A_j\prod_{j=1}^{k-1}B_jA_j^{-1}B_j^{-1}=1.
\label{mona-lisa}
\end{equation*}
Note that
\begin{eqnarray*}
B_j^\tau&=&B_j^{-1}\\
A_j^\tau&=&B_jA_jB_j^{-1}.
\end{eqnarray*}


\begin{thebibliography}{99}
\bibitem{abr}\ M. B. Abrahamse, R. G. Douglas. A class of subnormal operators
related to multiply connected domains. Adv. in Math., 19 (1976) 106--148. 

\bibitem{all}\ N. L. Alling. Extensions of meromorphic
function rings over non-compact Riemann surfaces I. Math. Z. 89 (1965) 273--299. 
 
\bibitem{vin}\ D. Alpay, V. Vinnikov.
Indefinite Hardy spaces on finite bordered Riemann surfaces.  
J. Funct. Anal. 172 (2000), no. 1, 221--248. 

\bibitem{aziz}\ T. Ya. Azizov, I. S. Iohvidov. Foundations of the
theory of linear operators in spaces with indefinite metric. Nauka,
Moscow, 1986, (in Russian). English translation: Linear operators in
spaces with an indefinite metric. John Wiley, New-York, 1989. 

\bibitem{ball}\ J. A. Ball, K. Clancey. Reproducing kernels for Hardy spaces
on multiply connected domains. Integral Equation Operator Theory. 25: 35--37,
1996. 

\bibitem{balvin}\ J. A. Ball, V. Vinnikov. Zero-pole interpolation for
meromorphic
matrix functions on an algebraic curve and transfer  functions of $2D$ systems.
Acta. Appl. Math., 45: 239--316, 1996. 

\bibitem{vvv4}\ J. A. Ball, V. Vinnikov.
Zero-pole interpolation for matrix meromorphic functions on a compact Riemann 
surface and a matrix Fay trisecant identity.  Amer. J. Math.  121  (1999),  
no. 4, 841--888. 

\bibitem{bali}\ J. A.  Ball, V. Vinnikov. 
Hardy spaces on a finite bordered Riemann surface, 
multivariable operator model theory and Fourier analysis along a unimodular curve. 
Systems, approximation, singular integral operators, and related topics 
(Bordeaux, 2000), 37--56, Oper. Theory Adv. Appl., 129, Birkhäuser, Basel, 2001.

\bibitem{bra}\ L. de Branges. Espaces hilbertiens de fonctions
 enti${\grave {\rm e}}$res. Masson, Paris, 1972. 

\bibitem{bogn}\ J. Bogn\'ar.
 Indefinite inner product spaces. Springer-Verlag,
Berlin, 1974. 


\bibitem{drez}\ J.-M. Drezet, M. S. Narasimhan. Groupe de Picard des
varetes de modules de fibres semi-stables sur les courbes algebriques.
 Invent. Math. 97 (1989), 53--94. 

\bibitem{dur}\ P. L. Duren. Theory of $H^p$ spaces. Vol. 38 of
 Pure and Applied
Math. Academic Press, New York, 1970. 

\bibitem{dym}\ H. Dym. J contractive matrix functions, reproducing
kernel spaces and interpolation. Vol. 71
 of CBMS Lecture Notes. Amer. Math. Soc.,
Rhodes Island, 1989. 

\bibitem{fkr}\ H. M. Farkas, I. Kra. Riemann surfaces. Second edition.
 Springer-Verlag. 1991. 


\bibitem{fay1}\ J. D. Fay. Theta functions on Riemann surfaces. Vol.
352 of Lecture Notes in Math. Springer-Verlag, New York, 1973. 

\bibitem{fay2}\ J. D. Fay. Kernel functions, analytic torsion,
and moduli spaces. Memoirs of the  Amer. Math. Soc. 464, Amer. Math.
Soc., Providence, 1992. 

\bibitem{forst}\ O. Forster. Lectures on Riemann surfaces. Springer-Verlag.
1991. 

\bibitem{garn}\ J. B. Garnett. Bounded analytic  functions. Academic press,
San Diego, 1981. 

\bibitem{grh}\ P. Griffits, J. Harris. Principles of algebraic geometry.
Wiley-Interscience Publication. 1994. 

\bibitem{gun1}\ R. C. Gunning. Lectures on Riemann surfaces.
Princeton mathematical notes. Princeton University press. 1966. 

\bibitem{hart}\ R. Hartshorne. Algebraic geometry. Springer-Verlag. 1977. 

\bibitem{hof}\ K. Hoffman. Banach spaces of analytic functions.
 Prentice-Hall, Englewood Cliffs, N.J. 1962. 

\bibitem{ioh}\ I. S. Iohvidov, M. G. Kre\u{\i}n, H. Langer. Introduction
to  the spectral theory  of operators in spaces with an indefinite metric.
 Akademie--Verlag, Berlin, 1982. 

\bibitem{kra}\ N. Kravitsky. Rational operator functions 
and Bezoutian operator vessels. Integral Equations Operator 
Theory 26 (1996), no. 1, 60--80.

\bibitem{mtz1}\ G. Mason, M. P. Tuite,  A. Zuevsky.  
Torus $n$-Point Functions for $\mathbb{R}$-graded Vertex Operator Superalgebras and Continuous
 Fermion Orbifolds.  Commun. Math. Phys.  vol 283 (2008) n.2, 305--342.

\bibitem{mtz11}\ M. P. Tuite,  A. Zuevsky. The Szeg\"o Kernel on a Sewn Riemann Surface. To appear. 

\bibitem{mtz2}\ M. P. Tuite,  A. Zuevsky. The genus two partition function for free fermionic vertex
operator algebras. To appear. 

\bibitem{mum}\ D. Mumford. Tata lectures on theta. Boston,
Birkh${\ddot {\rm a}}$user, 1982--1984. 

\bibitem{nar}\ R. Narasimhan. Compact Riemann surfaces.
Birkh${\ddot {\rm a}}$user. 1992. 

\bibitem{natan}\ S. M. Natanzon. Moduli spaces of real curves. Trans.
Moscow. Math. Soc., 37: 233--272, 1980. 

\bibitem{rudin}\ W. Rudin.  Analytic functions of class $H_p$. Trans.
Amer. Math. Soc., 78: 46--66, 1955. 

\bibitem{sarason}\ D. Sarason.
 The $H^p$ spaces of an annulus. Mem. Amer. Math. Soc., 1 (56),
1965. 

\bibitem{sesh}\ C. S. Seshadri. Fibr\'es
 vectoriels sur les courbes alg\'ebriques.
 Ast\'erisque 96. Socit\'et\'e  math\'ematique de  France. 1980. 

\bibitem{span}\ E. H.  Spanier. Algebraic topology. McGraw-Hill
book company. 1966. 

\bibitem{vvv1}\ V. Vinnikov. Complete
description of determinantal representations
of smooth irreducible curves. Linear Algebra Appl., 125: 103--140, 1989. 

\bibitem{vvv2}\ V. Vinnikov. Commuting nonselfadjoint operators and algebraic
 curves. In T. Ando and I. Gohberg, editors, Operator theory and complex analysis,
volume 59 of Operator Theory: Adv. Appl., 348--371.
Birkh${\ddot {\rm a}}$user Verlag, Basel, 1992. 


\bibitem{vvv5}\ V. Vinnikov. Self-adjoint determinantal representations of
real plane curves. Math. Ann. 296: 453--479, 1993. 

\bibitem{vvv3}\ V. Vinnikov. $2D$ systems and realization of bundle mappings
on compact Riemann surfaces. In U. Helmke, R. Mennicken  and J.Saurer, editors,
Systems and Networks: Mathematical theory  and applications (Vol II),
 volume 79 of
Math. Res., 909--912. Akademie Verlag, Berlin, 1994. 


\bibitem{vino}\ V. Vinnikov. Commuting operators and function
theory on a Riemann surface. Holomorphic spaces (S. Axler at al., eds),
Math. Sci. Res. Inst. Publ. 33, Cambridge Univ. Press, 1998. 

\bibitem{voi1}\ M. Voichik. Ideals  and invariant subspaces of analytic
functions. Trans. Amer. Math. Soc., 111: 493--512, 1964. 
\bibitem{y}\ A. Yamada. Precise variational formulas for abelian
differentials. Kodai Math. J. \textbf{3} (1980) 114--143.
\end{thebibliography}
\end{document}